\documentclass[3p,times,12pt]{elsarticle}
\usepackage{graphicx}
\usepackage{longtable}
\usepackage{amsmath,amssymb,amsthm,amsfonts}
\usepackage{xcolor}

\newcommand{\Spec}{\mathrm{Spec}}

\newtheorem{theorem}{Theorem}[section]
\newtheorem{prop}[theorem]{Proposition}
\newtheorem{lemma}[theorem]{Lemma}

\newtheorem{cor}[theorem]{Corollary}

\begin{document}

\begin{frontmatter}
\title{On the minimum number of distinct eigenvalues of a threshold graph}

\author{ Shaun Fallat, Seyed Ahmad Mojallal}
\address{Department of Mathematics and Statistics, University of Regina, Regina, Saskatchewan, S4S0A2,  Canada}
\address{ \,shaun.fallat@uregina.ca,\,ahmad\_mojalal@yahoo.com}
\begin{abstract} 
For a graph $G$, we associate a family of real symmetric matrices, $S(G)$, where for
any $A\in S(G)$, the location of the nonzero off-diagonal entries of $A$ are governed by
the adjacency structure of $G$. Let $q(G)$ be the minimum number of distinct eigenvalues over all matrices in $S(G)$. In this work, we give a characterization of all connected threshold graphs $G$ with $q(G)=2$. Moreover, we study the values of $q(G)$ for connected threshold graphs with trace $2$, $3$, $n-2$, $n-3$, where $n$ is the order of threshold graph.  The values of $q(G)$ are determined for all connected threshold graphs with $7$ and $8$ vertices with two exceptions. Finally, a sharp upper bound for $q(G)$ over all connected threshold graph $G$ is given.

\end{abstract}

\end{frontmatter}

\section{Introduction}

Let $G=(V(G), E(G))$ be a simple graph with vertex set $V(G)=\{1,2, \ldots, n\}$ and edge set $E(G)$. The order the graph $G$ is defined as $|V(G)|$. The set $S(G)$ of symmetric matrices described by $G$ consists of the set of all real symmetric $n \times n$ matrices $A=(a_{ij})$ such that for $i\neq j$, $a_{ij}\neq 0$ if and only if $ij\in E(G)$. We denote the spectrum of $A$, i.e., the multiset of eigenvalues of $A$, by $\Spec(A)$.  Given a graph $G$, the spectral invariant $q(G)$ is defined as follows:
 $$q(G)=\min\{q(A)\,:\, A\in S(G)\},$$ where $q(A)$ is the number of distinct eigenvalues of $A$.
The spectral invariant $q(G)$ is called the {\em minimum number of distinct eigenvalues of the graph $G$}. 

All matrices are real and symmetric; $O$ and $I$ denote zero and identity matrices of appropriate size, respectively. If the distinct eigenvalues of $A$ are $\lambda_1< \lambda_2< \ldots <\lambda_q$ and the multiplicities of these eigenvalues are $m_1, m_2, \ldots, m_q$ respectively, then the ordered multiplicity list of $A$ is $m(A) =(m_1, m_2, \ldots, m_q)$ and, additionally, the spectrum of $A$ is denoted by 
$$Spec(A)=\{\lambda_1^{[m_1]},\, \lambda_2^{[m_2]},\, \ldots,\,\lambda_q^{[m_q]}\}.$$ 
Here a simple eigenvalue $\lambda_k^{[1]}$ is simply denoted by  $\lambda_k$.
 The maximum multiplicity of $G$ is $M(G)=max\{mult_A(\lambda) :A\in S(G), \lambda\in \Spec(A)\}$, and the minimum rank of $G$ is $mr(G)=min\{rank\, A;\,A\in S(G)\}$. If we let $S_+(G)$ be all positive semidefinite matrices belongs to $S(G)$, then the maximum positive semidefinite nullity of $G$ is $M_+(G)=max\{mult_A(\lambda) :A\in S_+(G), \lambda\in \Spec(A)\}$. A subgraph $H$ of a graph $G$ is a graph whose vertex set and edge set are subsets of those of $G$. If $H$ is a subgraph of $G$, then $G$ is said to be a supergraph of $H$. 

The class of matrices $S(G)$ has been of interest to many researchers recently (see
\cite{Fal, Fer} and the references therein), and there has been considerable development of the
parameters $M(G)$  and $mr(G)$, and their positive semidefinite counterparts, see, for example, the
works \cite{Boo, FH, Fal}. Furthermore, as a consequence interest has grown in connecting these
parameters to various combinatorial properties for a given graph. For example, the inverse eigenvalue problem for graphs (see \cite{Hog}) continues to receive considerable and deserved
attention, as it remains one of the most interesting unresolved issues in combinatorial
matrix theory.

Threshold graphs can be characterized in many ways. One way of obtaining
a threshold graph is through an iterative process which starts with an isolated vertex, and where, at each step, either a new isolated vertex is added, or a vertex adjacent to all previous vertices (dominating vertex) is added. We may represent a threshold graph on $n$ vertices using a binary sequence $(b_1, \ldots, b_n)$. Here $b_i$ is 0 if vertex $v_i$ was added as an isolated vertex, and $b_i$ is 1 if $v_i$ was added as a dominating vertex. This representation has been called a {\em creation sequence} \cite{HSS}. For convenience, we use 0 as the first character of the string; it represents the first vertex of the graph. The number of characters equal to 1 in this string, called the {\em trace} of the graph, indicates the number of dominating vertices in its construction \cite{Mer}. As we are interested in connected threshold graphs, we insist $b_n =1$. In general, the trace $T$ of a graph $G$ is maximum index $i$ such that $d_i\ge i$, where $d_i$ is $i$th largest degree of $G$. Finally, for brevity, if $G$ is a threshold graph with creation sequence $(b_1, \ldots, b_n)$, then we write $G \cong (b_1, \ldots, b_n)$.

 For a threshold graph $G$, let $A\in S(G)$. We denote the submatrix of $A$ with rows indexed by dominating vertices and columns indexed by isolated vertices in $G$ by $A(\bar{1},\,\bar{0})$. Obviously, $A(\bar{1},\,\bar{0})$ is an $T\times (n-T)$ matrix. The submatrix of $A$ obtained by deleting the row and column corresponding to a vertex $i$ is denoted by $A(i)$ and the submatrix of $A$ obtained by deleting its $i$th row is denoted by $A(i;)$. 
 A split graph is a graph whose vertices can be partitioned into a clique and an independent set. In a split graph, if each vertex of the clique is adjacent to each vertex of the independent set, the graph is called a complete split graph.

Our main interest in this work lies in studying the value of $q(G)$ for connected threshold graphs. It is clear that knowledge of $q(G)$ for a graph $G$ will impact current studies of the inverse eigenvalue problem for graphs, and, in particular, the parameters $M(G)$ and $mr(G)$.

Recently, J.  Ahn et al. \cite{AAB} offered a complete solution to the ordered multiplicity inverse eigenvalue problem for graphs on six vertices. This work, represents in part, our motivation to study the values of $q(G)$ for connected threshold graphs with $7$ and $8$ vertices as well as several classes of threshold graphs with arbitrary orders.

The paper is organized as follows. In Section~\ref{Vertex-clique incidence matrix}, the notion of the vertex-clique incidence matrix of a graph $G$ is introduced and an application of this notion in the theory of minimum number of distinct eigenvalues of a graph is explored here. In Section~\ref{section:notation}, a list of well-known results and facts are given.
In Section~\ref{q=2}, a complete characterization of all threshold graphs with $q(G)=2$  is given. 
 In Section~\ref{specified threshold graphs}, the minimum number of distinct eigenvalues for a few specified threshold graphs are
 determined. 
   In Section~\ref{Threshold graphs with T=2}, a full characterization of $q(G)$ of all threshold graphs with trace two is given. In the next three sections we study $q(G)$ for threshold graphs with trace equal to $3$ and complete characterizations are provided for trace equal to $n-2$ and $n-3$.
  In Section~\ref{n=7} and \ref{n=8}, we calculate the values of $q(G)$ for any connected threshold graph of order $n=7$ and for all such threshold graphs of order 8, with two exceptions.
 Finally, 
 in Section~\ref{upper bound}, we provide  a sharp upper bound for $q(G)$ of any connected threshold graph of order $n$.

 \section{Vertex-clique incidence matrix of a graph}\label{Vertex-clique incidence matrix}

In this section, we introduce the notion of vertex-clique incidence matrix of a graph $G$ and its application to the issue of determining the minimum number of distinct eigenvalues of a graph. 

  An edge clique cover $F$ of a simple and undirected graph $G$ is a set of cliques of $G$, which together contain each edge of $G$ at least once.  The smallest cardinality of any edge clique cover of $G$ is called the edge clique cover number of $G$, and is denoted by $cc(G)$. 
  This number exists as the edge set of $G$ forms an edge clique cover for $G$. 
 An edge clique cover of $G$ with size $cc(G)$ is referred to as a  minimum edge clique cover of $G$.

 The vertex-clique incidence matrix of $G$ associated with the edge clique cover $F$ is defined as follows:
 Corresponding to any clique partition $F=\{C_1, C_2,\ldots, C_k\}$, we define a real matrix $M_F$ with rows and columns indexed by the vertices in $V$ and the cliques in $F$, respectively, such that the $ij$-entry of $M_F$ is real and nonzero if and only if the vertex $i$ belongs to the clique $C_j\in F$. 
 To ensure $M_F\,M^T_F\in S(G)$, we arrange that the inner product of row $i$ and column $j$ (when $i \neq j$) in $M_F$ is nonzero if and only if $ij\in E(G)$.

Recall that if $A$ and $B$ are any matrices of orders $t \times s$ and $s \times t$, respectively, then $AB$ and $BA$ have the same nonzero eigenvalues. Thus in order to minimize the number of distinct eigenvalues of $M_F\,M_F^T$, we need to minimize the number of distinct eigenvalues of $M_F^T\,M_F$. This is the key idea used to generalize the vertex-clique incidence matrix obtained from an edge clique cover by considering any positive real entries for $M_F$ or any negative real entries for $M_F$ but carefully preserving the condition that $M_F\,M_F^T\in S(G)$.

 \section{Notations and preliminaries}\label{section:notation}
 
In this section, we shall list some previously known results that are needed as well as some notation used throughout the paper.

Let $\circ$ denote the Schur (also known as the Hadamard or entrywise) product. 
The symmetric matrix $A$ has the {\it Strong Spectral Property} (or $A$ has the SSP for short) if the only symmetric matrix $X$ satisfying $A\circ X=O$, $I\circ X=O$ and $[A,\, X]=AX-XA=O$ is $X=O$.

The following result is given in \cite[Thm. 10]{BFH}.
\begin{lemma}{\rm \cite{BFH}}\label{ssp}
If $A\in S(G)$ has the SSP, then every supergraph of $G$ with the same vertex set has a realization that has the same spectrum as $A$ and has the SSP.
\end{lemma}

The following definitions and results can be found in \cite{LOS}:
For a graph $G$ and a vertex $v\in V(G)$, let  $N(v)$ denotes the open neighbourhood of $v$, i.e., the set of all vertices in $G$ adjacent to $v$. Let $dup(G,v)$ be a graph obtained by duplicating $v$, that is, augmenting $G$ with a new vertex $w$ and extra edges joining $w$ to every vertex in $N(v)$, so that in $dup(G,v)$ we have $N(v)=N(w)$ and $v$ and $w$ are not neighbours. Alternatively,  ``joined duplicated vertex" graph $jdup(G, v)$  is equal to $dup(G, v)$ with an extra edge joining $v$ and its duplicate vertex, $w$.  
\begin{lemma} {\rm \cite{LOS}}\label{jdup}
If $G$ is a non-empty graph, then $$q(jdup(G,\, v))\le q(G).$$
\end{lemma}

\begin{lemma}{\rm \cite{LOS}}\label{dup}
If $v$ is a vertex of $G$ and $A\in S(G)$ with $\lambda=a_{v,v}\in \sigma(A)$, then there is a matrix $C\in S(dup(G,v))$ with $\lambda=c_{v,v}$ so that the spectra of $A$ and $C$ are equal as sets, with the multiplicity of $\lambda$ in $\sigma(C)$ increased by 1, and so that $A$ is equal to $C$ with the $w$th row and column deleted, where $w$ is the duplicate vertex of $v$ introduced to form $dup(G,v)$ from $G$. 
\end{lemma}

Similar results can also be found in \cite{AAB}.

\begin{lemma}{\rm \cite{AAB}}\label{clon1}
Let $G$ be a graph with $M\in S(G)$ having multiplicity list $(\gamma_1,\, \ldots, \gamma_i,\ldots, \gamma_k)$ where the eigenvalue 0 has multiplicity $\gamma_i$. Then the following two cases are possible:\\
1.\, If the diagonal entry of $M$ corresponding to $v_j$ is zero, then the graph $H$ attained from $G$ by cloning $v_j$ without an edge has a matrix $N\in S(H)$ that attains the multiplicity list 
$(\gamma_1\, \ldots, \gamma_i+1,\, \ldots, \gamma_k)$. \\
2.\, If the diagonal entry of $M$ corresponding to $v_j$ is nonzero, then the graph $H$ attained from $G$ by cloning $v_j$ with an edge has a matrix $N\in S(H)$ that attains the multiplicity list 
$(\gamma_1\, \ldots, \gamma_i+1,\, \ldots, \gamma_k)$. 
\end{lemma}

Recall that the length of a path is simply the number of edges in that path, and that the distance between two vertices,
(in the same component) is the length of the shortest path between those two vertices. The next result represents a basic combinatorial lower bound on $q(G)$. 
\begin{lemma}{\rm \cite{AACF} }\label{diam}
 If there are vertices $u, \,v$ in a connected graph $G$ at distance $d$ and the path of length $d$ from $u$ to $v$ is unique, then $q(G) \ge d+1$.
\end{lemma}

Let $G$ be the star on $n$ vertices (that is, $G=K_{1,n-1}$, the complete bipartite graph with disjoint parts having sizes 1, n-1, respectively). Using Lemma, \ref{diam}, it is clear that $q(G) \geq 3$. On the other hand, the adjacency matrix associated with $G$ has the eigenvalues $\pm \sqrt{n-1}, 0$. Hence $q(G)=3$. Further observe that $G$ is a threshold graph with creation sequence $G \cong (0,0,\ldots, 0, 1)$, where there are $n-1$ consecutive zeros in this string. Hence for any threshold graph $H$ with creation sequence $(0,0,\ldots 0, 1,1, \ldots1)$, where the number of 0's and 1's is positive satisfies $q(H)\leq q(G) \leq 3$, by application of Lemma \ref{jdup}.

 \section{Threshold graphs with $q(G)=2$}\label{q=2}

In this section we provide a characterization for connected threshold graphs that satisfy $q(G)=2$.

\begin{theorem}\label{cosubmatrix}
For a threshold graph $G$, $q(G)=2$ if and only if there exists a matrix $A\in S(G)$ such that $A(\bar{1},\,\bar{0})$ is column orthogonal.
\end{theorem}
\begin{proof} It is well-known that if $q(G)=2$ then there exists an $n \times n$ orthogonal matrix $A\in S(G)$. This gives any $n \times k$ submatrices of $A$ must be column orthogonal. In particular, we consider the submatrix of $A$ with rows indexed by all vertices and columns indexed by only isolated vertices. This gives  $A(\bar{1},\,\bar{0})$ is column orthogonal.

Now, assume that there exists a matrix $A\in S(G)$ s.t. $A(\bar{1},\,\bar{0})$ is column orthogonal. This leads to a vertex-clique incidence matrix $M$ of the form $\left(
  \begin{array}{c}
    A(\bar{1},\,\bar{0}) \\
    X\\
  \end{array}
\right),$ where $X=diag(x_1, x_2, \ldots, x_k)$ and $k=n-T$. We can choose the values of $x_i$ under which the size of each column is constant $C$. This gives $M^T\,M=C^2\,I_k$, and then 
$q(MM^T)=2$. To complete the proof, it suffices to show that $M\,M^T\in S(G)$.    

Assume that $v$ is a dominating vertex and $w$ an isolated vertex. We have $vw\in E(G)$ if and only if the inner product between rows corresponding to $v$ and $w$ in $M$ is nonzero. If both $v$ and $w$ are isolated, then the inner product of their rows in $M$ are zero, that is, $x_ie_i^T\,.\, x_je_j^T=0$. Now, we show that if $v$ and $w$ are dominating vertices in $G$, then inner product between their rows in $M$ is always nonzero. Suppose this is not true and there exist some pairs of dominating vertices s.t.  their inner products are equal to zero. Now considering a proper diagonal matrix $D=diag(d_1, d_2, \ldots, d_k)$, we can make a new matrix $B=A(\bar{1},\,\bar{0})\,D$ s.t.  $B^T\,B$ is a $k\times k$ diagonal matrix and $B\,B^T$ is a $T\times T$ matrix with all nonzero entries. Hence we arrive at the new vertex-clique incidence matrix $\overline{M}=\left(
  \begin{array}{c}
    B \\
    Y\\
  \end{array}
\right),$ where $Y=diag(y_1, y_2, \ldots, y_k)$ s.t. $\overline{M}\,\overline{M}^T\in S(G)$ and $\overline{M}^T\,\overline{M}=aI_k$, for a positive $a$.  This completes the proof of the theorem. 
 \end{proof}

We now apply Theorem \ref{cosubmatrix} to certain families of threshold graphs.

\begin{theorem}\label{compare1}
Let $G\cong (\underbrace{0, \ldots, 0}_{k_1},\underbrace{1, \ldots, 1}_{t_1}, \underbrace{0, \ldots, 0}_{k_2},\underbrace{1, \ldots, 1}_{t_2}, \ldots,  \underbrace{0, \ldots, 0}_{k_s}, \underbrace{1, \ldots, 1}_{t_s})$ be a connected threshold graph with $n$ vertices and let $k_i, t_i, s\ge 1$. If $k_i\ge \sum_{j=i}^s t_j$ for an integer $i\in \{2,\ldots,s\}$, then $q(G)\ge 3$. 
\end{theorem}
\begin{proof} Since threshold graph $G$ is connected, we have $q(G)\ge 2$. Suppose $q(G)=2$. Then there exists $A\in S(G)$ s.t. $A^2=I$. Now, let $B$ is a submatrix of $A$ of order $n\times (1+k_i)$, with rows indexed by all vertices and columns indexed by $1+k_i$ isolated vertices, the first isolated vertex together with $k_i$ related isolated vertices.  Obviously, since $A$ is orthogonal, $B$ is a column orthogonal matrix and then the $t\times (1+k_i)$ submatrix $B^{\prime}$ of $B$, where $t=\sum_{j=i}^s t_j$. By our assumption, we have $1+k_i>t$, that is, the matrix $B^{\prime}$ is a $t\times (1+k_i)$ with the number of columns greater than the number of rows, a contradiction since $B^{\prime}$ is column  orthogonal. Hence $q(G)\neq 2$. That is, $q(G)\ge 3$.     
 \end{proof}

 \begin{theorem}\label{}
Let $G$ be a connected threshold graph of order $n$ and trace $T$. If $q(G)=2$, then $T\ge \lceil \frac{n}{2} \rceil$.
\end{theorem}
\begin{proof}
Since $q(G)=2$, there exists $A \in S_{+}(G)$ such that $m(A)=(a,b)$ with 
 $\max\{a,\, b\}\le T $. This gives $n=a+b\le 2T$, and thus $T\ge \lceil \frac{n}{2} \rceil$.
\end{proof}
This leads to the following result:
\begin{cor}\label{trace}
Let $G$ be a connected threshold graph of order $n$ and trace $T$. If $T<\lceil \frac{n}{2} \rceil$, then $q(G)\ge 3$. 
\end{cor}

\section{On $q(G)$ for specific families of threshold graphs}\label{specified threshold graphs}

In this section we determine $q(G)$ for threshold graphs of the form $G\cong (0,1,0,1, \ldots,0,1)$, \\ $G\cong (0,\underbrace{1, \ldots, 1}_{t_1},0,\underbrace{1, \ldots, 1}_{t_2},\ldots,0,\underbrace{1, \ldots, 1}_{t_s})$ and for  the special class of complete split graphs. 

\begin{theorem}\label{thm00}
 Let $G\cong (0,1,0,1, \ldots,0,1)$ be a threshold graph of order $n$. Then $q(G)=3$. Moreover, for any supergraph $H$ of $G$ on the same vertex set satisfies $q(H)\le 3$.
\end{theorem}
\begin{proof} We consider the set $C$ of cliques in $G$ with each clique corresponding to each zero in (or isolated vertex in $G$) in the creation sequence for $G$. Clearly, 
$C$ is a minimum edge clique cover for $G$. Now, we form the vertex-clique incidence matrix  $M_1$ of $G$ corresponding to $C$: $$M_1=\left(
  \begin{array}{ccccccc}
    1 & 0 & 0 & \ldots & 0  \\
    1 & 0 & 0 &  \ldots & 0  \\
    1 & 2 & 0 &  \ldots & 0   \\
 1 & 1 & 3 &  \ldots & 0   \\
   \vdots & \vdots & \ddots &  \ddots   & \vdots \\
    1 & 1 &  \ldots& 1 & t  \\
    0 & 1 & 0 & \ldots & 0  \\
    0& 0 & 1 & \ddots & \vdots \\
    \vdots & \vdots &\ddots  & \ddots & 0  \\
    0 & 0 &  \ldots & 0 & 1 \\
  \end{array}
\right),$$ where $t=\frac{n}{2}$ is the number of zeros in the creation sequence for $G$. It is not difficult to see that the square of the length of the last column of $M_1$ is $t^2+1$. We can arrange that the squares of the lengths of all other columns also equal $t^2+1$ by setting a new element given by $x_i$ given in $i$th column in $M_2$ for $1\le i \le t-1$ as follows:
$$M_2=\left(
  \begin{array}{ccccccc}
    x_1 & 0 & 0 & \ldots & 0  \\
    1 & 0 & 0 &  \ldots & 0  \\
    1 & 2 & 0 &  \ldots & 0   \\
 1 & 1 & 3 &  \ldots & 0   \\
   \vdots & \vdots & \ddots &  \ddots   & \vdots \\
    1 & 1 &  \ldots& 1 & t  \\
    0 & x_2 & 0 & \ldots & 0  \\
    0& 0 & x_3 & \ldots & 0 \\
    \vdots & \vdots &\ddots  & \ddots & \vdots  \\
    0 & 0 & 0 &  \ldots & 1 \\
  \end{array}
\right).$$ The matrix of $M_2^T M_2$ is a $t \times t$ matrix with each entry on the main diagonal equal to $t^2+1$  and all off-diagonal entries equal to $t$. Thus the spectrum of $M_2^T M_2$ is $$\Spec(M_2^T M_2)=\left\{t^2-t+1^{[t-1]},\,2t^2-t+1\right\}.$$
Then the spectrum of $A=M_2M_2^T$ is equal to $$\Spec(A)=\left\{t^2-t+1^{[t-1]},\,2t^2-t+1, 0^{[t]}\right\}.$$
Obviously, $A\in S(G)$, and hence $q(G)\le 3$. Now, using Lemma \ref{diam} by considering the unique distance $d$ between the unique pendant vertex and any vertex other than the unique universal vertex, we arrive at $q(G)\ge 3$. Thus $q(G)=3$ and the first part of the proof is complete.

It is not difficult to determine (using mathematical software, for example) that $A$ has SSP and thus by Lemma \ref{ssp}, any supergraph $H$ of $G$ on the same vertex set has a realization with the same spectrum as $A$, which gives $q(H)\le 3$. 
\end{proof}

\begin{theorem}\label{salter}
Let $G\cong (0,\underbrace{1, \ldots, 1}_{t_1},0,\underbrace{1, \ldots, 1}_{t_2},\ldots,0,\underbrace{1, \ldots, 1}_{t_s})$, where $t_i\ge 1$ for $i = 1,2,\ldots, s$, and $s\ge 2$. If $t_s=1$, then $q(G)=3$; otherwise if $t_s\ge 2$, then $q(G)=2$.
\end{theorem}
\begin{proof} First we assume that $t_s=1$. By Lemma \ref{diam}, $q(G)\ge 3$. On the other hand, by Lemma \ref{jdup} and Theorem \ref{thm00} we have $q(G)\le q(H)=3$, where $H \cong (0,1,0,1,\ldots,0,1)$. Hence $q(G)=3$.

Next, we assume that $t_s\ge 2$.  Let $H\cong (0,1,0,1,\ldots,0,1,0,1,1)$. A vertex-clique incidence matrix of $H$ is as follows:
$$M=\left(
  \begin{array}{cccccc}
    1 & 0 & 0 & \ldots & 0 & 0 \\
    1 & -s+1 & 0 & \ddots & 0 & 0 \\
    1 & 1 & -s+2 & 0 & \ddots & 0 \\
    \vdots & \vdots  & \vdots  & \ddots &\ddots   & 0 \\
    1 & 1 & 1 & 1 & \ldots & -1\\
    1 & 1 & 1 & 1 & \ldots & 1\\
    1 & 0 & 0 & 0& \ldots & 0 \\
    0 & 1 & 0 & \ddots & \ddots & 0 \\
    0 & 0 & 1 & 0 & \ddots & \vdots \\
    \vdots & \vdots  & \vdots  & \ddots &\ddots   & 0 \\
    0 & 0 & 0 &  \ldots & 1 & 0 \\
    0 & 0 & 0 &  \ldots & 0 & 1 \\
  \end{array}
\right).$$ We have $A=M M^T\in S(H)$ and $M^T M=(s^2-s+1)\,I_{s}$. Thus $\Spec(A)=\Big\{0^{s+1},\,{s^2-s+1}^{[s]}\Big\}$. Hence, Lemma \ref{jdup} $q(G)\le q(H)=2$ and consequently, $q(G)=2$.
 \end{proof}

In the following we obtain the minimum number of distinct eigenvalues for complete split graphs.

\begin{theorem}\label{csplit}
Let $G\cong (\underbrace{0, \ldots, 0}_{k_1},\underbrace{1, \ldots, 1}_{t_1})$, where $t_1, k_1\ge 1$. If $k_1\le t_1$, then $q(G)=2$; otherwise $q(G)=3$.
\end{theorem}
\begin{proof}  By Theorem \ref{cosubmatrix}, $q(G)=2$ if and only if there exists a matrix $A\in S(G)$ such that $A(\bar{1},\,\bar{0})$ is column orthogonal. For the threshold graph $G$, the matrix $A(\bar{1},\,\bar{0})$ has the following form:
$$A(\bar{1},\,\bar{0})=\left(
  \begin{array}{cccccc}
    a_1 & b_1 &  \ldots & c_1 \\
    a_2 & b_2 &  \ddots &  c_2\\
    \vdots & \vdots  & \ddots &\vdots \\
    a_{t_1} & b_{t_1}  & \ldots &  c_{t_1}\\
  \end{array}
\right),$$ consisting of all nonzero entries. Such a matrix $A(\bar{1},\,\bar{0})$ allows a column orthogonal matrix if $k_1\le t_1$. Otherwise, $q(G)\ge 3$. This with Lemma \ref{jdup} gives $q(G)=3$. 
 \end{proof}

\section{Threshold graphs with trace $T=2$}\label{Threshold graphs with T=2}

In this section, we completely characterize $q(G)$ for all connected threshold graphs $G$ with trace equal to two.  Note that any such graph must have at least 3 vertices. If $G$ is a threshold graph with trace equal to two and has 3 vertices, then $G=K_3$, and $q(G)=2$. If such a graph has 4 vertices, then either $G \cong (0,1,0,1)$ or $G \cong (0,0,1,1)$. In the later case $G$ is isomorphic to the complete graph on 4 vertices with one edge removed, and it is well-known that $q(G)=2$ in this case. For the former case, it follows that $q(G)=3$ from \cite{AACF}. For most of the results in this section, we assume that $n \geq 5$. 

We begin with a series of facts by considering specific trace two threshold graphs. 

\begin{theorem}\label{thm060}
Let $G\cong (\underbrace{0, \ldots, 0}_{k},1,0,1)$, where $k\ge 1$. Then $q(G)=3$.
\end{theorem}
\begin{proof} We assume that $H\cong (0,1,0,1)$ with a vertex set as $(1,2,4,3)$ with the presented order respecting the order in the creation sequence. We consider the following clique-incidence matrix of $H$:$$M=\left(
  \begin{array}{cccc}
\sqrt{3} & 0 \\
    1 &  0  \\
    1 &  2 \\
    0 & 1 \\
  \end{array}
\right).$$ Then $M\,M^T\in S(H)$. On the other hand, 
$M^T\,M=5\,I_2+2\,A(K_2)$, where $A(K_2)$ is the adjacency matrix of the complete graph of order 2. Thus $\Spec(M^T\,M)=\{7,\, 3\}$, and consequently, $\Spec(M\,M^T)=\{7,\, 3,\, 0^{[2]}\}$. The eigenvalue $3$ is also an entry on the main diagonal of $M\,M^T$ corresponding to vertex 1. Then applying Lemma \ref{dup} by duplicating vertex 1, we arrive at $q(G)\le 3$. This combined with Lemma \ref{diam} completes the proof of theorem.  
\end{proof}

\begin{theorem}\label{thm061}
Let $G\cong (0,1,\underbrace{0, \ldots, 0}_{k},1)$, where $k\ge 1$. Then $q(G)=3$.
\end{theorem}
\begin{proof} We assume that $H\cong (0,1,0,1)$ with a vertex set as $(1,2,4,3)$ with the presented order respecting the order in the creation sequence. We consider the following clique-incidence matrix of $H$:$$M=
\left(  \begin{array}{cccc}
  1 &  0 \\
    1 &  0 \\
    1 &  1 \\
    0 & \sqrt{2} \\
  \end{array}
\right).$$ Then $M\,M^T\in S(H)$. On the other hand,
$M^T\,M=3\,I_2+A(K_2)$, where $A(K_2)$ is the adjacency matrix of the complete graph of order 2. Thus $\Spec(M^T\,M)=\{4,\, 2\}$, and consequently, $\Spec(M\,M^T)=\{4,\, 2,\, 0^{[2]}\}$. The eigenvalue $2$ is also an entry on the main diagonal of $M\,M^T$ corresponding to the vertex 4. Then applying Lemma \ref{dup} by duplicating  vertex 4, we arrive at $q(G)\le 3$. This combined with Lemma \ref{diam} completes the proof of theorem. 
\end{proof}

\begin{theorem}\label{thm01}
 Let $G\cong (0,0,1, \underbrace{0, \ldots, 0}_{k}, 1)$ be a threshold graph of order $n$ with $k\ge 1$. Then $q(G)=3$. 
\end{theorem}
\begin{proof} Let  $H\cong (0,0,1,0,1)$ be a threshold graph on five vertices. A vertex-clique incidence matrix $M$ of $H$ is as follows:
$$M=\left(
  \begin{array}{cccc}
    3 & 0 & 0 \\
    1 & -3 & 0  \\
   1 & 1 & -2 \\
    0 & 1 &0 \\
    0 & 0 & \sqrt{7} \\
 \end{array}
\right).$$ The matrix $A=MM^T\in S(H)$, and the spectrum of $A$ is
$$\Spec(A)=\{13^{[2]},\,0^{[2]},\,7\}.$$ Then $q(H)\le 3$. This with Lemma \ref{diam} gives $q(H)=3$. On the main diagonal of $A$, the entry corresponding to the unique pendant vertex $v$ is equal to $7$ and since $7\in \Spec(A)$. Thus any threshold graph $G$ obtained from $H$ with duplicating the pendant vertex $v$ as often as is required, has a realization with the same spectrum as $A$ by Lemma \ref{dup}. Then $q(G)\le 3$. Again using Lemma \ref{diam}, we conclude $q(G)=3$.
\end{proof}

\begin{theorem}\label{thm10}
 Let $G\cong (\underbrace{0, \ldots, 0}_{k_1}, \underbrace{1, \ldots, 1}_{t_1}, \underbrace{0, \ldots, 0}_{k_2}, \underbrace{1, \ldots, 1}_{t_2})$ be a threshold graph of order $n$ with $k_1,\, k_2, t_1, t_2 \ge 1$. Then $q(G)\le 4$. 
 \end{theorem}
\begin{proof} Let  $H\cong (0,1,0,0,1)$ be a threshold graph on five vertices. A vertex-clique incidence matrix $M$ of $H$ is as follows:
$$M=\left(
  \begin{array}{cccc}
    1 & 0 & 0 \\
    1 & 0 & 0  \\
    1 & 1 & 1 \\
    0 & 1 &0 \\
    0 & 0 & 1 \\  
  \end{array}
\right).$$ The matrix $A=MM^T\in S(H)$. The spectrum of $A$ is
$$\Spec(A)=\{4.4142,\,1.5857,\,1,\,0^{[2]}\}.$$ Then $q(H)\le 4$. On the main diagonal of $A$, the entries corresponding to the first and forth vertices in  $H\cong (0,1,0,0,1)$ are equal to $1$ and since $1\in \Spec(A)$, any threshold graph $G' \cong (\underbrace{0, \ldots, 0}_{k_1}, 1, \underbrace{0, \ldots, 0}_{k_2}, 1)$ obtained from $H$ by duplicating these vertices as needed, has a realization with the same spectrum as $A$ using Lemma \ref{dup}. Then $q(G')\le 4$. Finally, applying Lemma \ref{jdup}, it follows that $q(G) \leq q(G') \leq 4$. Hence the proof is complete.
\end{proof}

Before we state and prove our main result in this section we preface that discussion by alerting the reader that the  threshold graph shown in Figure 1 is the only threshold graph on 7 vertices which attains $q(G)=4$ (see Section 10), all other threshold graphs $H$ on 7 vertices satisfy $q(H)\leq 3$.

\begin{center}
 \includegraphics[height=5cm,keepaspectratio]{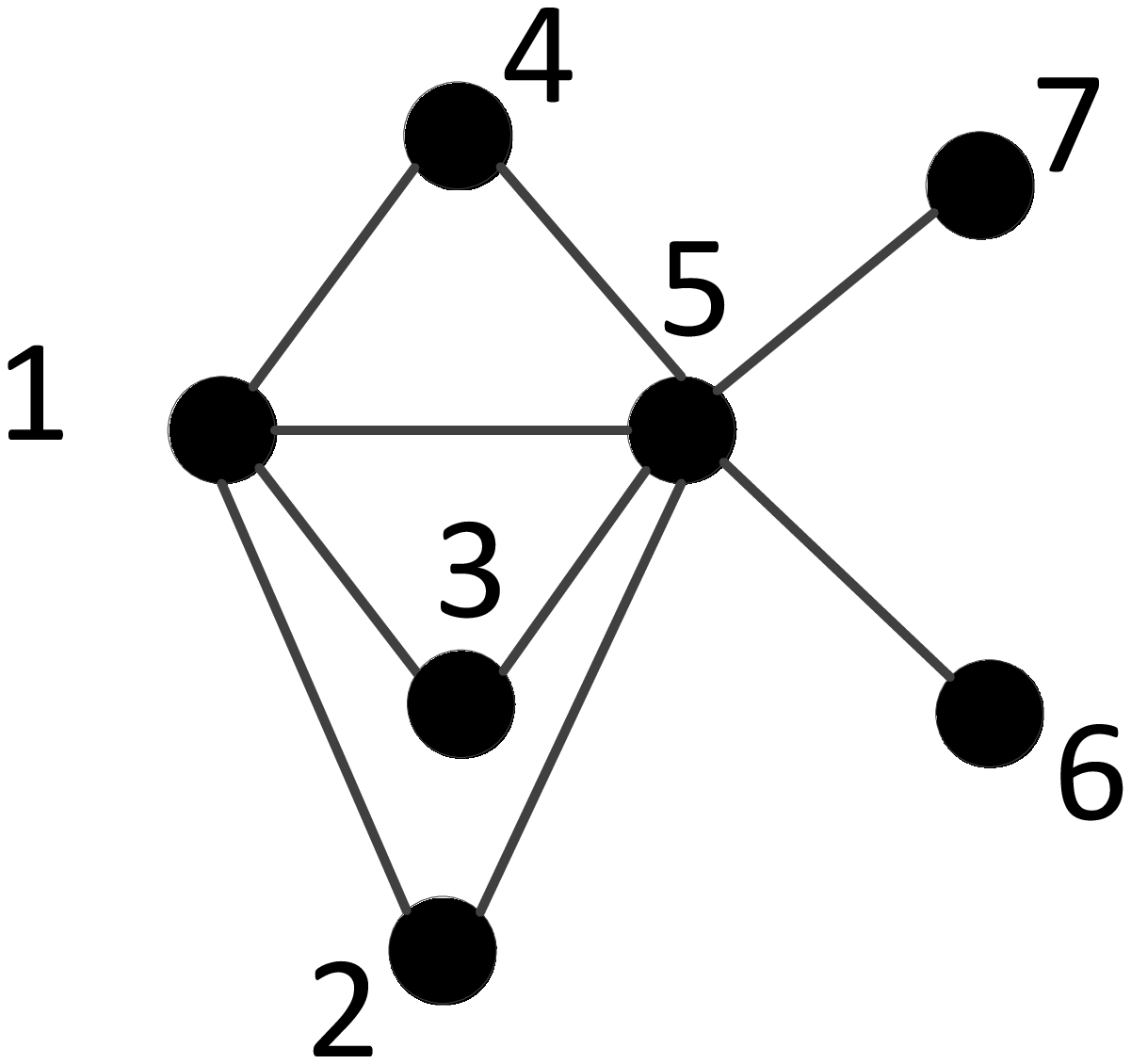}\\
  $G$\\
 \vspace*{1cm}
 \noindent
 Fig.1. The threshold graph $G$.
 \end{center}

In the next result we state and prove the main theorem of this section.
\begin{theorem}\label{T=2}
 Let $G\cong (\underbrace{0,\ldots,0}_{k_1},1, \underbrace{0,\ldots,0}_{k_2},1)$ be a threshold graph of order $n\ge 5$ with $k_1\ge 1$ and $k_2\ge 0$. Then
$$q(G)=\begin{cases}
3~~~~\mbox{if}\, k_1\in \{1,\,2\} ~\mbox{or}~ k_2\in \{0,\,1\},\\
4~~~~\mbox{otherwise}.
\end{cases}$$
 \end{theorem}
\begin{proof}
Considering the values of the parameters $k_1$ and $k_2$, we may conclude that $q(G)=3$ whenever $k_1\in \{1,\,2\} ~\mbox{or}~ k_2\in \{0,\,1\}$ by applying Theorems \ref{csplit}, \ref{thm060}, \ref{thm061} and \ref{thm01}. So we may assume that the parameters $k_1$ and $k_2$ fall outside of these ranges. Noting, Theorem \ref{thm10}, we have that any such threshold graph $G$ satisfies $q(G) \leq 4$. 

The next step in the proof involves the particular trace 2 threshold graph $G \cong (\underbrace{0,\ldots,0}_{k-1},1,0,0,1)$, with $k \geq 4$. We will show that $q(G)=4$.
By Lemma \ref{diam} we have $q(G)\ge 3$. We show that the equality is not possible. If $q(G)=3$ then since we know that largest possible positive semidefinite nullity is 2, since the trace of $G$ is 2. Furthermore, from \cite{ShuWan} we know that the largest possible multiplicity is $n-4$. Thus it follows that the only possible allowed ordered multiplicity list for $A\in S(G)$ with $q(A)=3$ is $(2, n-4,2)$. This follows since the multiplicities of the largest and smallest eigenvalues is at most 2, as the clique cover number of $G$ is the trace of $G$, which equals 2.
Without loss of generality  we may assume such that $A=(a_{i,j})_{1\le i,j\le n}\in S(G)$ has spectrum $$\Spec(A)=\{-\lambda^{[2]},\,0^{[n-4]},\, 1^{[2]}\},~~~~\mbox{where}~ 0 <\lambda < 1.$$ 
The main diagonal entries $a_{i,i}$ must be equal to zero for $i=2,\ldots,k$, otherwise the nullity of zero will be changed to $n-3$ by Lemma \ref{clon1}, which is a contradiction with the fact $M(G)=n-4$ \cite{ShuWan}. For convenience, we write $A$ in the following form: 
$$A=\left(
  \begin{array}{ccccccc}
    d_0 & c_1 & \cdots & c_{k-1} & c_k & 0 & 0 \\
    c_1 & 0   & \cdots &  0 & c_{k+1} & 0 & 0\\
    \vdots & 0  & \ddots & 0 & \vdots& \vdots  & \vdots\\
    c_{k-1} & 0 & \cdots & 0&  c_{2k-1} & 0 & 0\\
    c_k & c_{k+1}&  \cdots & c_{2k-1} &  d_k & c_{2k}& c_{2k+1} \\
    0  & 0 &  \cdots & 0&  c_{2k} & d_{k+1} & 0 \\
    0  & 0 &  \cdots  & 0&  c_{2k+1} & 0 & d_{k+2} \\
  \end{array}
\right),$$ where all $c_j$ are necessarily non-zero, $k=n-3$, and the last 3 vertices correspond to the last three positions in the creation sequence for $G$ in reverse order.
Since $rank(A)=4$, then either $d_{k+1}=d_{k+2}=0$ or both $d_{k+1}$ and $d_{k+2}$ are non-zero. Consequently, we consider these two cases separately.

{\bf Case\,1)}\, $d_{k+1}=d_{k+2}=0$. First we show that $$-\lambda=d_0-1,~~\mbox{and}~~~d_0=d_{k}.$$ To do this, we assume that $B=A(k+1)$. The spectrum of $B$ is
$$\Spec(B)=\left\{\frac{1}{2}d_0\pm \frac{1}{2}\sqrt{4\sum_{i=1}^{k-1}c_i^2+d_0^2},\, 0^{[n-4]}\right\}.$$ Applying interlacing for $A$ and $B$ we obtain
 \begin{align}
 \frac{1}{2}d_0+\frac{1}{2}\sqrt{4\sum_{i=1}^{k-1}c_i^2+d_0^2} &=1 \nonumber\\
                                                               & \label{equA10}\\
 \frac{1}{2}d_0-\frac{1}{2} \sqrt{4\sum_{i=1}^{k-1}c_i^2+d_0^2} &=-\lambda.\nonumber
 \end{align}
This proves $-\lambda=d_0-1$ and $\sum_{i=1}^{k-1}c_i^2=1-d_0$. Considering the trace of $A$ we have
$d_0+d_k=2(d_0-1)+2$ and then $d_0=d_k$. 

Again applying interlacing for $A$ and $C=A(1)$, we arrive at 
\begin{equation}\label{equ20}
\sum_{i=k+1}^{2k+1} c_i^2=1-d_0.
\end{equation} Comparing the coefficients of $x^{n-2}$ in the characteristic polynomial of $A$ and $(x-d_0+1)^2\, (x-1)^2\, x^{n-4}$
we conclude that 
$$-\sum_{i=1}^{2k+1} c_i^2+d_0\,d_k=d_0^2+2d_0-2.$$ This equation combined with the fact $d_0=d_k$ and (\ref{equ20}) gives $$-[2(\sum_{i=1}^{k-1}c_i^2)+c_k^2]=2d_0-2,$$ which implies $c_k=0$, a contradiction. \\

{\bf Case\,2)}\, Both $d_{k+1}$ and $d_{k+2}$ are non-zero. Following as in the beginning of Case 1 above, we may conclude that 
\begin{align}
 \frac{1}{2}d_0+\frac{1}{2}\sqrt{4\sum_{i=1}^{k-1}c_i^2+d_0^2} &=1 \nonumber\\
                                                               & \label{equA10}\\
 \frac{1}{2}d_0-\frac{1}{2} \sqrt{4\sum_{i=1}^{k-1}c_i^2+d_0^2} &=-\lambda.\nonumber
 \end{align} 
 This follows since if the smallest or largest eigenvalue of a symmetric matrix coincides with a main diagonal entry, then the remainder of the row and column containing that diagonal entry must be zero. Furthermore it now follows from interlacing applied for $A$ and $A(k+1)$, that $d_{k+1}$ and $d_{k+2}$ have different signs. We may assume that $d_{k+1}<0<d_{k+2}$.
This proves $-\lambda=d_0-1$ and $\sum_{i=1}^{k-1}c_i^2=1-d_0$. Furthermore, we claim the following conditions hold for the matrix $A$: $c_{i+k}=bc_i$, for some $b \neq 0$ and $i=1,2,\ldots, k-1$; $d_k = \frac{c_{2k}^2}{d_{k+1}} + \frac{c_{2k+1}^2}{d_{k+2}} +2bc_k-b^2d_0$, and finally $c_k=b(d_0-1)$; The first two conditions follow by applying elementary row/column operations to $A$ and noting that the rank of $A$ is 4 and that both $d_{k+1}$ and $d_{k+2}$ are non-zero. To verify the third condition, since $-\lambda$ and 1 are multiple eigenvalues of $A$, we have $-\lambda$ and 1 are both eigenvalues of $A(k+1)$. Since $A(k+1)$ is a direct sum of a matrix $C$ and a $2 \times 2$ diagonal matrix consisting of main diagonal entries $d_{k+1}$ and $d_{k+2}$, it follows that $-\lambda$ and 1 are both eigenvalues of $C$. Applying row and column operations to the matrix $A(k+2, k+3)-I$, we may conclude that $c_k=b(d_0-1)$.

Label the two pendant vertices in $G$ as $x$ and $y$. We let $A(x)$ (and $A(y)$) denote the principal submatrix of $A$ obtained by deleting the row and column corresponding to $x$ ($y$). Further, let $B = A(x,y)$. Under the assumptions above it follows that the rank of $B$ is at most 3, and if we let $B(v)$ denote the principal submatrix of $B$ obtained by deleting the last row and column of $B$, then the rank of $B(v)$ is 2. Now define $f_{x}(z) = \det(zI-A(x))$, and $f_{y}(z) = \det(zI-A(y))$. Applying basic properties of determinants we have 
\[ f_{x}(z) = (z-d_{k+2})f_{B}(z) - c_{2k+1}^2f_{B(v)}(z), \;\;\;\; f_{y}(z) = (z-d_{k+1})f_{B}(z) - c_{2k}^2f_{B(v)}(z),\]
where $f_B$ and $f_{B(v)}$ are the characteristic polynomials of $B$ and $B(v)$ respectively. Since $-\lambda$ is a multiple eigenvalue of $A$ we have $f_x(-\lambda)=f_y(-\lambda)=0$. Similarly, since neither 
$d_{k+1}$ nor $d_{k+2}$ can equal $-\lambda$ it follows that $-\lambda$ is an eigenvalue of $B(v)$. Thus
\[ 0 = f_x(-\lambda)-f_{y}(-\lambda) = f_B(-\lambda) (d_{k+1} - d_{k+2}). \] As noted above that rank of $B$ is at most 3 and so the characteristic polynomial of $B$ has the form:

\[ f_{B}(z) = z^{k+1} - (d_0+d_k)z^{k}+\alpha z^{k-1} + \beta z^{k-2}, \]
where
\[ \alpha = -\sum_{i \neq k}c_i^2 + d_0d_k-c_k^2, \;\;\; \beta = \sum_{i=1}^{k-1} (2c_ic_kc_{k+i}-d_kc_i^2 - d_0c_{k+1}^2).\]
Using the established conditions on the entries of $A$ above, we may simplify both $\alpha$ and $\beta$ as follows:
\[\alpha =-(1-d_0 + b^2(1-d_0)) -d_0d_k -b^2(d_0-1)^2, \;\;\; \beta = (1-d_0)(b^2(2(d_0-1)-d_0)-d_k).\]
Using these simplified equations, we may now write:
\[ f_B(-\lambda) (d_{k+1} - d_{k+2}) = (d_{k+1}-d_{k+2})(d_0-1)^{k-1}(b)^2(d_0-2)^2.\] 
A contradiction arises since the latter expression cannot equal 0. Thus no such $A \in S(G)$ with $q(A)=3$ exists, so $q(G)=4$. 

More generally, assume that $G\cong (\underbrace{0,\ldots,0}_{k_1},1, \underbrace{0,\ldots,0}_{k_2},1)$, where $k_2>2$. Suppose there exists $A \in S(G)$ with $q(A)=3$. Then, as before we may assume that the rank of $A$ is 4 and that the spectrum of $A$ is $$\Spec(A)=\{-\lambda^{[2]},\,0^{[n-4]},\, 1^{[2]}\},~~~~\mbox{where}~ 0 <\lambda < 1.$$ Now suppose all of the diagonal entries corresponding to the pendant vertices of $G$ are zero in $A$. Then, as in Case 1 above, we arrive at a contradiction. On the other hand, if at least one main diagonal entry corresponding to a pendant vertex is nonzero, then there are at least two such nonzero main diagonal entries. However, in this case it follows that the rank of $A$ is at least 5 since the number of such pendant vertices is at least 3. This completes the proof.
\end{proof}

\begin{cor}\label{pa}
Let $G\cong (0,\underbrace{1, \ldots, 1}_{t_1},\underbrace{0, \ldots, 0}_{k_2},1)$, where $t_1,\,k_2\ge 1$ or let 
$G\cong (0,1,\underbrace{0, \ldots, 0}_{k_2},\underbrace{1, \ldots, 1}_{t_2})$, where $t_t,\,k_2\ge 1$ 
Then in either case we have $q(G)=3$.
\end{cor}

\section{Threshold graphs with trace $T=3$}\label{T=3}

In this section, we study $q(G)$ for connected threshold graphs with trace equal to three.
\begin{theorem}\label{thm05}
Let $G\cong (\underbrace{0, \ldots, 0}_{k_1},1,0,1,\underbrace{0, \ldots, 0}_{k_2},1)$, where $k_1,\,k_2\ge 1$. Then $q(G)=3$.
\end{theorem}
\begin{proof} We assume that $H\cong (0,1,0,1,0,1)$ with a vertex set labelled as $(1,2,5,3,6,4)$ respecting the order of vertices in its creation sequence. We consider the following clique-incidence matrix of $H$:$$M=\left(
  \begin{array}{cccc}
\sqrt{7} & 0 & 0 \\
    \sqrt{2} & 0 & 0  \\
    1 & -3 & 0 \\
    1 & 1 &-2   \\
    0 & 1 &  0 \\
    0 & 0 & \sqrt{7} \\
  \end{array}
\right).$$ Observe that  $M\,M^T\in S(H)$, and note that
$M^T\,M=11\,I_3-2\,A(K_3)$, where $A(K_3)$ is the adjacency matrix of the complete graph of order 3. Hence $\Spec(M^T\,M)=\{13^{[2]},\, 7\}$, and consequently, $\Spec(M\,M^T)=\{13^{[2]},\, 7,\, 0^{[3]}\}$. The eigenvalue $7$ is also an entry on the main diagonal of $M\,M^T$ corresponding to the vertices 1 and 6. Applying Lemma \ref{dup} by duplicating the vertices 1 and 6, we arrive at $q(G)\le 3$. This combined with Lemma \ref{diam} completes the proof of theorem. 
\end{proof}

\begin{theorem}\label{thm06}
Let $G\cong (0,1,\underbrace{0, \ldots, 0}_{k},1,0,1)$, where $k\ge 1$. Then $q(G)=3$.
\end{theorem}
\begin{proof} We assume that $H\cong (0,1,0,1,0,1)$ with a vertex set labelled as $(1,2,5,3,6,4)$ respecting the order of vertices in its creation sequence. We consider the following clique-incidence matrix of $H$:$$M=\left(
  \begin{array}{cccc}
\sqrt{2} & 0 & 0 \\
    1 & 0 & 0  \\
    1 & 1 & 0 \\
    1 & 1 & 2   \\
    0 & \sqrt{3} &  0 \\
    0 & 0 & 1 \\
  \end{array}
\right).$$ 
The remainder of the proof is very similar to the previous proof and is omitted.
\end{proof}


\begin{theorem}\label{q<4}
 Let $G\cong (\underbrace{0, \ldots, 0}_{k_1}, 1, \underbrace{0, \ldots, 0}_{k_2}, 1, \underbrace{0, \ldots, 0}_{k_3}, 1)$ be a threshold graph of order $n$ with $k_1,\,k_2,\,k_3\ge 1$. Then $q(G)\le 4$. 
\end{theorem}
\begin{proof} Let  $H\cong (0,1,0,1,0,0,1)$ be a threshold graph on seven vertices. A vertex-clique incidence matrix $M$ associated with $H$ is as follows:
$$M=\left(
  \begin{array}{cccc}
    1 & 0 & 0 & 0 \\
    1 & 0 & 0 & 0 \\
    1 & 1 & 0 & 0 \\
    1 & 1 & 1 & 1 \\
    0 & \sqrt{2} & 0 & 0 \\  
    0 & 0 & 1 & 0 \\  
    0 & 0 & 0 & 1 \\  
  \end{array}
\right).$$ The matrix $A=MM^T\in S(H)$, and the spectrum of  $A$
$$\Spec(A)=\{7,\,2^{[2]},\,0^{[3]},\,1\}.$$ Thus $q(H)\le 4$. On the main diagonal of $A$, the entries corresponding to the first, sixth and seventh vertices in  $H\cong (0,1,0,1,0,0,1)$ are equal to $1$ and $2$ and since $\{1,\,2\} \subseteq \Spec(A)$, any threshold graph $G$ obtained from $H$ by duplicating these vertices as often as needed, has a realization with the same spectrum as $A$ by Lemma \ref{dup}. Then $q(G)\le 4$, which completes the proof.
\end{proof}

\section{Threshold graphs with trace $T=|G|-2$}\label{Threshold graphs with T=n-2}

In this section, we calculate $q(G)$ for all threshold graphs $G$ with trace equal to $n-2$, where $n$ is the number of vertices of $G$.

\begin{theorem}\label{T=n-2}
Let $G\cong (0,\underbrace{1, \ldots, 1}_{t_1},0,\underbrace{1, \ldots, 1}_{t_2})$, where $t_1\ge 0$ and $t_2\ge 1$. 
\begin{enumerate}
\item If $t_1=0$ and $t_2\ge 2$ then $q(G)=2$;
\item  if $t_1=0$ and $t_2=1$, then $q(G)=3$;
\item if $t_1\ge 1$ and $t_2=1$, then $q(G)=3$; and 
\item if $t_1\ge 1$ and $t_2\ge 2$, then $q(G)=2$.
\end{enumerate}
\end{theorem}
\begin{proof} The result for the case $t_1=0$, is directly confirmed by Theorem \ref{csplit}. Now, suppose $t_1\ge 1$. First assume that $t_2=1$. Then by Lemma \ref{diam} , we have $q(G)\ge 1+d=3$. On the other hand, by Lemma \ref{jdup}, $q(G)\le q(0,1,0,1)=3$, which proves the desired result in this case.  Second, assume that $t_2\ge 2$. Consider the graph $H\cong (0,1,0,1,1)$ with the following vertex-clique incidence matrix:
$$M=\left(
  \begin{array}{cc}
     1 & 0 \\
     2 & 1 \\
    -2 & 1 \\
     1 & 0 \\
     0 & \sqrt{8} \\
  \end{array}
\right).$$ Then $A=M\,M^T\in S(H)$ and $M^T\,M=10\,I_2$. Hence $\Spec(A)=\{0^{[3]},10^{[2]}\}$ and consequently, $q(H)=2$. By Lemma \ref{jdup} we conclude that $q(G)\le q(H)=2$ and thus $q(G)=2$. 
 \end{proof}

\section{Threshold graphs with trace $T=|G|-3$}\label{Threshold graphs with T=n-3}

In this section, we compute the values of $q(G)$ for all connected threshold graphs $G$ with  trace equal to $n-3$, where $n$ is the number of  vertices of $G$.

\begin{theorem}\label{T=n-3}
Let $G$ be a connected threshold graph with $n$ vertices and trace $T=n-3$. 
\begin{enumerate}
\item For $$G\cong (0,\underbrace{1, \ldots, 1}_{t_1},0,\underbrace{1, \ldots, 1}_{t_2},0,\underbrace{1, \ldots, 1}_{t_3}),$$  we have 
$q(G)=3$ if $t_3=1$ and $q(G)=2$, otherwise;
\item  for $G\cong (0,0,0,\underbrace{1, \ldots, 1}_{t_1})$, $q(G)=2$ if $t_1\ge 3$, and $q(G)=3$ if $t_1=1,2$; 
\item for $G\cong (0,0,\underbrace{1, \ldots, 1}_{t_1},0,\underbrace{1, \ldots, 1}_{t_2})$, $q(G)=3$ if $t_2=1$ and $q(G)=2$, if $t_2\ge 2$; and 
\item for $G\cong (0,\underbrace{1, \ldots, 1}_{t_1},0,0,\underbrace{1, \ldots, 1}_{t_2})$, $q(G)=3$ if $t_2=1,2$ and $q(G)=2$ if $t_2\ge 3$. 
\end{enumerate}
\end{theorem}
\begin{proof} For the threshold graph $G\cong (0,\underbrace{1, \ldots, 1}_{t_1},0,\underbrace{1, \ldots, 1}_{t_2},0,\underbrace{1, \ldots, 1}_{t_3})$ with $t_i\ge 1$, the desired conclusions follow from Theorem \ref{salter}.

For the threshold graph $G\cong (0,0,0,\underbrace{1, \ldots, 1}_{t_1})$, the conclusions follow from Theorem \ref{csplit}.

Now let $G\cong (0,0,\underbrace{1, \ldots, 1}_{t_1},0,\underbrace{1, \ldots, 1}_{t_2})$ with $t_2=1$. By Lemma \ref{diam} , we have $q(G)\ge 1+d=3$. On the other hand, by Lemma \ref{jdup} and Theorem \ref{T=2}, $q(G)\le q(0,0,1,0,1)=3$. This gives $q(G)=3$. Otherwise, $t_2\ge 2$. 
Now, consider the graph $H\cong (0,0,1,0,1,1)$ with the following vertex-clique incidence matrix:
$$M=\left(
  \begin{array}{ccc}
     1 & -2 & 0 \\
     1 & 1 & -1 \\
     1 & 1 & 1 \\
     2 & 0 & 0 \\
     0 & 1 & 0 \\
     0 & 0 & \sqrt{5}\\
  \end{array}
\right).$$ Then $A=M\,M^T\in S(H)$ and $M^T\,M=7\,I_3$. Hence $\Spec(A)=\{0^{[3]},7^{[3]}\}$ and consequently, $q(H)=2$. Now, by Lemma \ref{jdup} we conclude that $q(G)\le q(H)=2$ or $q(G)=2$.

Finally, let $G\cong (0,\underbrace{1, \ldots, 1}_{t_1},0,0,\underbrace{1, \ldots, 1}_{t_2}).$ If  $t_2=1$, then by Lemma \ref{diam} , we have $q(G)\ge 3$. On the other hand, by Lemma \ref{jdup} and Theorem \ref{T=2}, $q(G)\le q(0,1,0,0,1)=3$. This gives $q(G)=3$. If $t_2=2$, then 
 by Theorem \ref{cosubmatrix}, $q(G)=2$ if and only if there exists a matrix $A\in S(G)$ such that $A(\bar{1},\,\bar{0})$ is column orthogonal. For the threshold graph $G$, the matrix $A(\bar{1},\,\bar{0})$ has the following form:
$$A(\bar{1},\,\bar{0})=\left(
  \begin{array}{cccccc}
    a_1 & 0 & 0 \\
    \vdots & \vdots  &  \vdots \\
    a_{t_1} & 0  & 0 \\
     a_{t_1+1} & b_1 & c_1 \\
     a_{t_1+2} & b_2 & c_2 \\
  \end{array}
\right),$$ where all variables above are required to be nonzero. Obviously, $A(\bar{1},\,\bar{0})$ allows a column orthogonal matrix if and only if $\left(
  \begin{array}{cccccc}
   a_{t_1+1} & b_1 & c_1 \\
   a_{t_1+2} & b_2 & c_2 \\
  \end{array}
\right)$ is column orthogonal, which is impossible as all entries are required to be nonzero. Then $q(G)\ge 3$. This with Lemma \ref{diam} and Theorem \ref{T=2} gives $q(G)=3$.

If $t_2\ge 3$, we consider the graph $H\cong (0,1,0,0,1,1,1)$. Again $q(H)=2$ if and only if there exists a matrix $A\in S(H)$ such that  $A(\bar{1},\,\bar{0})$ is column orthogonal.

For the threshold graph $H$, the matrix $A(\bar{1},\,\bar{0})$ has the following form:
$$A(\bar{1},\,\bar{0})=\left(
  \begin{array}{cccccc}
    a_1 & 0 & 0 \\
    a_2 & b_2  & c_2 \\
    a_3 & b_3  & c_3 \\
    a_4 & b_4 & c_4 \\
  \end{array}
\right),$$ where all variables above are required to be nonzero. 
It is easy to check that $A(\bar{1},\,\bar{0})$ allows a column orthogonal matrix if and only if $\left(
  \begin{array}{ccc}
  a_2 & b_2  & c_2 \\
    a_3 & b_3  & c_3 \\
    a_4 & b_4 & c_4 \\
  \end{array}
\right)$ is column orthogonal which is possible. Then $q(H)=2$. This with Lemma \ref{jdup}  gives $q(G)\le q(H)=2$, that is, $q(G)=2$.
 \end{proof}

\section{Threshold graphs with $7$ vertices}\label{n=7}

In this section, we determine the values of $q(G)$ for all connected threshold graphs $G$ with $7$ vertices.

 The list of  all connected threshold graphs $G$ of order $7$ along with their value $q(G)$ is given below.\\
\begin{tabular}{|c|c|c|c|c|c|} 
\hline
Index & Creation Sequence &  $q(G)$ & Explanations\\
\hline
1& (0,\,0,\,0,\,0,\,0,\,0,\,1) & 3  &  Star graph\\
\hline 
2& (0,\,1,\,0,\,0,\,0,\,0,\,1) &   3  &   Thm \ref{T=2}\\
\hline
3 & (0,\,0,\,1,\,0,\,0,\,0,\,1) &  3  &  Thm \ref{T=2}\\
\hline
4& (0,\,0,\,0,\,1,\,0,\,0,\,1) &   4  &   Thm \ref{T=2} \\
\hline
5& (0,\,0,\,0,\,0,\,1,\,0,\,1) &   3  &   Thm \ref{T=2} \\
\hline
6& (0,\,0,\,0,\,0,\,0,\,1,\,1) &    3  &   Thm \ref{T=2} \\
\hline
7& (0,\,1,\,1,\,0,\,0,\,0,\,1) &   3  &   Corollary \ref{pa}\\
\hline
8& (0,\,1,\,0,\,1,\,0,\,0,\,1) &   3  &   Thm \ref{thm05}\\
\hline
9& (0,\,1,\,0,\,0,\,1,\,0,\,1) &   3  &   Thm \ref{thm06} \\
\hline
10& (0,\,1,\,0,\,0,\,0,\,1,\,1) &   3  &  Lemma \ref{jdup}  \& Corollaries \ref{trace} \&  \ref{trace}\\
\hline
11& (0,\,0,\,1,\,1,\,0,\,0,\,1) &  3  &  Thm \ref{T=2} \& Lemmas \ref{diam} \& \ref{jdup}  \\
\hline
12& (0,\,0,\,1,\,0,\,1,\,0,\,1) &  3  &  Thm \ref{thm05}  \\
\hline
13&(0,\,0,\,1,\,0,\,0,\,1,\,1) &   3  &  Lemma \ref{jdup}  \& Corollary \ref{trace} \\
\hline
14&(0,\,0,\,0,\,1,\,1,\,0,\,1) &   3  &  Thm \ref{T=2} \& Lemmas \ref{diam} \& \ref{jdup}  \\
\hline
15&(0,\,0,\,0,\,1,\,0,\,1,\,1) &   3  &  Thm \ref{T=2} \& Lemma \ref{jdup}  \& Corollary \ref{trace}
\\
\hline
16&(0,\,0,\,0,\,0,\,1,\,1,\,1) &   3  &  Lemma  \ref{jdup}  \& Corollary \ref{trace} \\
\hline
17&(0,\,1,\,1,\,1,\,0,\,0,\,1) &   3  &  Thm \ref{T=2} \& Lemma \ref{jdup}  \\
\hline
18&(0,\,1,\,1,\,0,\,1,\,0,\,1) &   3  &  Thm \ref{thm00} \& Lemma \ref{jdup}  \\
\hline
19&(0,\,1,\,1,\,0,\,0,\,1,\,1) &   3  &  Thm \ref{T=n-3}  \\
\hline
20&(0,\,1,\,0,\,1,\,1,\,0,\,1) &   3  & Thm \ref{thm00} \& Lemma \ref{jdup}  \\
\hline
21&(0,\,1,\,0,\,1,\,0,\,1,\,1) &    2  &  Thm \ref{T=n-3} \\
\hline
22& (0,\,1,\,0,\,0,\,1,\,1,\,1) &   2  &  Thm \ref{T=n-3}  \\
\hline
23& (0,\,0,\,1,\,1,\,1,\,0,\,1) &   3  &  Thm \ref{T=2} \& Lemma \ref{jdup}  \\
\hline
24& (0,\,0,\,1,\,1,\,0,\,1,\,1) &   2 &  Thm  \ref{T=n-3}  \\
\hline
25& (0,\,0,\,1,\,0,\,1,\,1,\,1) &   2  &  Thm  \ref{T=n-3}  \\
\hline
26& (0,\,0,\,0,\,1,\,1,\,1,\,1) &   2 & Thm  \ref{csplit}  \\
\hline
27& (0,\,1,\,0,\,1,\,1,\,1,\,1) &   2  & Thm \ref{T=n-2}  \\
\hline
28& (0,\,1,\,1,\,0,\,1,\,1,\,1) &   2  & Thm \ref{T=n-2}  \\
\hline
29& (0,\,1,\,1,\,1,\,0,\,1,\,1) &   2  &   Thm \ref{T=n-2} \\
\hline
30& (0,\,1,\,1,\,1,\,1,\,0,\,1) &   3  & Lemmas \ref{diam} \& \ref{jdup}  \\
\hline
31& (0,\,1,\,1,\,1,\,1,\,1,\,1) &   2  & Complete graph  \\
\hline
\end{tabular}

\section{Threshold graphs with $8$ vertices}\label{n=8}

In this section, we determine the values of $q(G)$ for all connected threshold graphs $G$ on $8$ vertices with two exceptions. For this, we need the following additional results.

\begin{prop}\label{a15a19}
Let $G_1\cong (0,\,0,\,1,\,0,\,0,\,1,\,0,\,1)$ and $G_2\cong (0,\,0,\,0,\,1,\,0,\,0,\,1,\,1)$.  Then $q(G_1)=q(G_2)=3$.
\end{prop}
\begin{proof}
Suppose that $H$ is the graph with the following vertex-clique incidence matrix:
$$M=\left(
  \begin{array}{ccccc}
    \sqrt{8} & 0 & 0 & 0 & 0 \\
     1 & -\sqrt{2} & 0 & 0 & 0\\
   \sqrt{2} & 1 & 0 & 0 & 0 \\
    0 & 1 & 1 & 1 & 0 \\
    0 & 0 & \sqrt{3} & 0 & 0 \\
    0 & 0 & 0 & \sqrt{3} & 0 \\
    0 & 1 & 1 & 1 & 2 \\
   0 & 0 & 0 & 0 & 1 \\
  \end{array}
\right).$$ Then $A=M\,M^T\in S(H)$. The matrix $A$ has SSP with the spectrum of  $\Spec(A)=\{0^{[3]},\,3^{[3]},\,11^{[2]}\}$ and consequently, $q(H)=3$. The threshold graphs $G_1$ and $G_2$ are supergraphs of $H$ with the same vertex set, which implies $q(G_i)\le 3$ for $i=1,2$, by Lemma \ref{ssp}. This with Corollary \ref{trace} gives the required result.  \end{proof}

\begin{prop}\label{31}
Let $G\cong (0,\,1,\,0,\,0,\,1,\,0,\,1,\,1)$. Then $q(G)=3$.
\end{prop}
\begin{proof} Let $A\in S(G)$ and $A(\bar{1},\,\bar{0})$ be the submatrix of $A$ with rows and columns corresponding to dominating and isolated vertices in $G$, respectively. Then this submatrix has  the following form:
$$A(\bar{1},\,\bar{0})=\left(
  \begin{array}{cccc}
     c_1 & 0 & 0 & 0 \\
     c_2 & c_3 & c_4 & 0 \\
    c_5 & c_6 & c_7 & c_8 \\
    c_9 & c_{10} & c_{11} & c_{12} \\
  \end{array}
\right)$$ with all nonzero entries $c_j$. It can be easily deduced that  $A(\bar{1},\,\bar{0})$ does not allow a column orthogonal matrix, so $q(G)\ge 3$ by Theorem \ref{cosubmatrix}. Now apply Theorem \ref{thm06} and Lemma \ref{jdup}, the desired result is obtained.
 \end{proof}

\begin{prop}\label{37}
Let $G\cong (0,\,0,\,1,\,0,\,1,\,0,\,1,\,1)$. Then $q(G)=2$.
\end{prop}
\begin{proof} Let $A\in S(G)$ and $A(\bar{1},\,\bar{0})$ be the submatrix of $A$ with rows and columns corresponding to dominating and isolated vertices in $G$, respectively. Consider
$$A(\bar{1},\,\bar{0})=\left(
  \begin{array}{cccc}
     1 & -3 & 0 & 0 \\
     1 & 1 & -2 & 0 \\
    1 & 1 & 1 & -1 \\
    1 & 1 & 1 & 1 \\
  \end{array}
\right).$$ Since $A(\bar{1},\,\bar{0})$ is  column orthogonal, by Theorem \ref{cosubmatrix}, $q(G)=2$. 
 \end{proof}

\begin{prop}\label{38}
Let $G\cong (0,\,0,\,1,\,0,\,0,\,1,\,1,\,1)$. Then $q(G)=2$.
\end{prop}
\begin{proof} Let $A\in S(G)$ and $A(\bar{1},\,\bar{0})$ be the submatrix of $A$ with rows and columns corresponding to dominating and isolated vertices in $G$, respectively. Then this submatrix must have  the following form:
$$B=A(\bar{1},\,\bar{0})=\left(
  \begin{array}{cccc}
     c_0 & c_1 & 0 & 0 \\
     c_2 & c_3 & c_4 & c_5 \\
    c_6 & c_7 & c_8 & c_9 \\
    c_{10} & c_{11} & c_{12} & c_{13} \\
  \end{array}
\right),$$ where all $c_i$ are required to be nonzero. It can be easily determined that  the $3\times 4$ submatrix $B(1;)$ obtained from $B$ with removal of its first row, allows a column orthogonal matrix. Then considering appropriate values for $c_0$ and $c_1$ under which $c_0\,c_1=-(c_2c_3+c_6c_7+c_{10} c_{11})$ we arrive at the column orthogonal matrix $B$. This with Theorem \ref{cosubmatrix} gives $q(G)=2$.
 \end{proof}

\begin{prop}\label{40}
Let $G\cong (0,\,0,\,0,\,1,\,1,\,0,\,1,\,1)$. Then $q(G)=2$.
\end{prop}
\begin{proof} Let $A\in S(G)$ and $A(\bar{1},\,\bar{0})$ be the submatrix of $A$ with rows and columns corresponding to dominating and isolated vertices in $G$, respectively. Then consider
$$A(\bar{1},\,\bar{0})=\left(
  \begin{array}{cccc}
     1 & 3/2 & 9/2 & 0 \\
     2 & 1/2 & -7/2 & 0 \\
    1/2 & -1 & 1 & 2 \\
    -1 & 2 & -2 & 1 \\
  \end{array}
\right).$$ Since $A(\bar{1},\,\bar{0})$ is column orthogonal, by Theorem \ref{cosubmatrix}, it follows that $q(G)=2$. 
 \end{proof}

 \begin{prop}\label{41}
Let $G\cong (0,\,0,\,0,\,1,\,0,\,1,\,1,\,1)$. Then $q(G)=2$.
\end{prop}
\begin{proof} 
 Suppose that $H$ is the graph given in Proposition \ref{40}. Using $A(\bar{1},\,\bar{0})$ used in the proof of Proposition \ref{40}, we can create the following vertex-clique incidence matrix for $H$:
$$M=\left(
  \begin{array}{ccccc}
     1 & 3/2 & 9/2 & 0  \\
     2 & 1/2 & -7/2 & 0 \\
     1/2 & -1 & 1 & 2 \\
    -1 & 2 & -2 & 1  \\
    1/2\sqrt{129} & 0 & 0 & 0  \\
    0 & \sqrt{31}  & 0 & 0 \\
    0 & 0 & 1  & 0 \\
   0 & 0 & 0 & \sqrt{67/2}  \\
  \end{array}
\right).$$ Then $A=M\,M^T\in S(H)$. The matrix $A$ has SSP with the spectrum of  $\Spec(A)=\{{\frac{77}{2}}^{[4]},\,0^{[4]}\}$ and consequently, $q(H)=2$. The threshold graph $G$ is a supergraph of $H$ and then $q(G)\le 2$, by Lemma \ref{ssp}. This with the fact that for any connected graph $G$,  $q(G)\ge 2$ gives the required result.
 \end{proof}

 The list of connected threshold graphs of order $8$ with associated values of $q(G)$ are listed below with two exceptions. \\
{\tiny \begin{tabular}{|c|c|c|c|c|c|} 
\hline
Index & Creation Sequence &  $q(G)$ & Explanations\\
\hline
1& (0,\,0,\,0,\,0,\,0,\,0,\,0,\,1) &     3  &  Star graph\\
\hline 
2& (0,\,1,\,0,\,0,\,0,\,0,\,0,\,1) &     3  &   Corollary \ref{pa}\\
\hline
3 & (0,\,0,\,1,\,0,\,0,\,0,\,0,\,1) &    3  &  Thm \ref{T=2}\\
\hline
4& (0,\,0,\,0,\,1,\,0,\,0,\,0,\,1) &     4  &    Thm \ref{T=2} \\
\hline
5& (0,\,0,\,0,\,0,\,1,\,0,\,0,\,1) &     4  &    Thm \ref{T=2} \\
\hline
6& (0,\,0,\,0,\,0,\,0,\,1,\,0,\,1) &     3  &    Thm \ref{T=2} \\
\hline
7& (0,\,0,\,0,\,0,\,0,\,0,\,1,\,1) &     3  &   Lemma \ref{jdup}  \& Corollary \ref{trace}\\
\hline
8& (0,\,1,\,1,\,0,\,0,\,0,\,0,\,1) &     3  &   Corollary \ref{pa}\\
\hline
9& (0,\,1,\,0,\,1,\,0,\,0,\,0,\,1) &     3  &   Thm \ref{thm05} \\
\hline
10& (0,\,1,\,0,\,0,\,1,\,0,\,0,\,1) &    3 or 4  &  Corollary \ref{trace} \& Thm \ref{q<4} \\
\hline
11& (0,\,1,\,0,\,0,\,0,\,1,\,0,\,1) &    3  &  Thm \ref{thm06}  \\
\hline
12& (0,\,1,\,0,\,0,\,0,\,0,\,1,\,1) &    3  &  Corollaries \ref{pa} \& \ref{trace}  \& Lemma \ref{jdup}  \\
\hline
13&(0,\,0,\,1,\,1,\,0,\,0,\,0,\,1) &    3  &  Lemma \ref{diam}  \& Corollary \ref{trace} \&  Thm \ref{T=2} \\
\hline
14&(0,\,0,\,1,\,0,\,1,\,0,\,0,\,1) &    3  &  Thm \ref{thm05}  \\
\hline
15&(0,\,0,\,1,\,0,\,0,\,1,\,0,\,1) &    3   &  Prop \ref{a15a19}
\\
\hline
16&(0,\,0,\,1,\,0,\,0,\,0,\,1,\,1) &    3  &  Lemma \ref{jdup}  \& Corollary \ref{trace} \&  Thm \ref{T=2} \\
\hline
17&(0,\,0,\,0,\,1,\,1,\,0,\,0,\,1) &    3 or 4  &  Corollary \ref{trace} \& Thm \ref{q<4} \\
\hline
18&(0,\,0,\,0,\,1,\,0,\,1,\,0,\,1) &    3   &  Thm \ref{thm05}  \\
\hline
19&(0,\,0,\,0,\,1,\,0,\,0,\,1,\,1) &    3   &  Prop \ref{a15a19} \\
\hline
20&(0,\,0,\,0,\,0,\,1,\,1,\,0,\,1) &    3  &  Corollary \ref{trace} \&  Thm \ref{T=2} \\
\hline
21&(0,\,0,\,0,\,0,\,1,\,0,\,1,\,1) &    3  &  Corollary \ref{trace} \&  Thm \ref{T=2} \\
\hline
22& (0,\,0,\,0,\,0,\,0,\,1,\,1,\,1) &   3  &  Lemma \ref{jdup}  \& Corollary \ref{trace}  \\
\hline
23& (0,\,1,\,1,\,1,\,0,\,0,\,0,\,1) &   3  &  Lemmas \ref{diam} \& \ref{jdup} \&  Thm \ref{T=2} \\
\hline
24& (0,\,1,\,1,\,0,\,1,\,0,\,0,\,1) &   3 &  Lemmas \ref{diam} \& \ref{jdup} \& Thm \ref{thm05}   \\
\hline
25& (0,\,1,\,1,\,0,\,0,\,1,\,0,\,1) &  3  &  Lemmas \ref{diam} \& \ref{jdup} \& Thm \ref{thm06}    \\
\hline
26& (0,\,1,\,1,\,0,\,0,\,0,\,1,\,1) &  3 &   Lemma \ref{jdup} \&  Thms \ref{T=2} \& \ref{compare1}  \\
\hline
27& (0,\,1,\,0,\,1,\,1,\,0,\,0,\,1) &   3  & Lemmas \ref{diam} \& \ref{jdup} \& Thm \ref{thm05}    \\
\hline
28& (0,\,1,\,0,\,1,\,0,\,1,\,0,\,1) &  3  & Thm \ref{thm00}  \\
\hline
29& (0,\,1,\,0,\,1,\,0,\,0,\,1,\,1) &   3 &  Lemma \ref{diam}  \& Thms \ref{thm05} \& \ref{compare1}   \\
\hline
30& (0,\,1,\,0,\,0,\,1,\,1,\,0,\,1) &   3  &  Lemmas \ref{diam} \& \ref{jdup} \& Thm \ref{thm06}  \\
\hline
31& (0,\,1,\,0,\,0,\,1,\,0,\,1,\,1) &   3  & Prop \ref{31} \\
\hline
32& (0,\,1,\,0,\,0,\,0,\,1,\,1,\,1) &   3  &  Lemma \ref{jdup} \&  Thms \ref{T=2} \& \ref{compare1}  \\
\hline
33& (0,\,0,\,1,\,1,\,1,\,0,\,0,\,1) &   3  & Lemmas \ref{diam} \& \ref{jdup} \& Thm \ref{T=2} \\
\hline
34& (0,\,0,\,1,\,1,\,0,\,1,\,0,\,1) &   3  &  Lemmas \ref{diam} \& \ref{jdup} \& Thm \ref{thm05}    \\
\hline
35& (0,\,0,\,1,\,1,\,0,\,0,\,1,\,1) &   3  & Lemma  \ref{jdup} \&  Thms \ref{T=2} \& \ref{compare1} \\
\hline
36& (0,\,0,\,1,\,0,\,1,\,1,\,0,\,1) &    3  &  Lemma \ref{diam}  \& Thm \ref{thm05}  \\
\hline
37& (0,\,0,\,1,\,0,\,1,\,0,\,1,\,1) &    2  &  Prop \ref{37}   \\
\hline
38& (0,\,0,\,1,\,0,\,0,\,1,\,1,\,1) &    2  &  Prop \ref{38} \\
\hline
39& (0,\,0,\,0,\,1,\,1,\,1,\,0,\,1) &    3  &  Lemmas \ref{diam} \& \ref{jdup} \&  Thm \ref{T=2} \\
\hline
40& (0,\,0,\,0,\,1,\,1,\,0,\,1,\,1) &    2 &  Prop \ref{40}\\
\hline
41& (0,\,0,\,0,\,1,\,0,\,1,\,1,\,1) &    2  &  Prop \ref{41} \\
\hline
42& (0,\,0,\,0,\,0,\,1,\,1,\,1,\,1) &    3  &  Thm  \ref{csplit}  \\
\hline
43& (0,\,1,\,1,\,1,\,1,\,0,\,0,\,1) &    3  & Thm \ref{T=n-3} \\
\hline
44& (0,\,1,\,1,\,1,\,0,\,1,\,0,\,1) &    3  & Thm \ref{T=n-3} \\
\hline
45& (0,\,1,\,1,\,1,\,0,\,0,\,1,\,1) &    3  & Thm \ref{T=n-3} \\
\hline
46& (0,\,1,\,1,\,0,\,1,\,1,\,0,\,1) &    3  &  Thm \ref{T=n-3} \\
\hline
47& (0,\,1,\,1,\,0,\,1,\,0,\,1,\,1) &    2  & Thm \ref{T=n-3} \\
\hline
48& (0,\,1,\,1,\,0,\,0,\,1,\,1,\,1) &   2  & Thm \ref{T=n-3} \\
\hline
49& (0,\,1,\,0,\,1,\,1,\,1,\,0,\,1) &   3  & Thm \ref{T=n-3} \\
\hline
50& (0,\,1,\,0,\,1,\,1,\,0,\,1,\,1) &   2  & Thm \ref{T=n-3} \\
\hline
51& (0,\,1,\,0,\,1,\,0,\,1,\,1,\,1) &   2  & Thm \ref{T=n-3} \\
\hline
52& (0,\,1,\,0,\,0,\,1,\,1,\,1,\,1) &   2  & Thm \ref{T=n-3} \\
\hline
53& (0,\,0,\,1,\,1,\,1,\,1,\,0,\,1) &   3  & Thm \ref{T=n-3} \\
\hline
54& (0,\,0,\,1,\,1,\,1,\,0,\,1,\,1) &   2  & Thm \ref{T=n-3} \\
\hline
55& (0,\,0,\,1,\,1,\,0,\,1,\,1,\,1) &   2  & Thm \ref{T=n-3} \\
\hline
56& (0,\,0,\,1,\,0,\,1,\,1,\,1,\,1) &   2  & Thm \ref{T=n-3} \\
\hline
57& (0,\,0,\,0,\,1,\,1,\,1,\,1,\,1) &   2  & Thm \ref{T=n-3} \\
\hline
58& (0,\,1,\,1,\,1,\,1,\,1,\,0,\,1) &   3  & Thm \ref{T=n-2} \\
\hline
59& (0,\,1,\,1,\,1,\,1,\,0,\,1,\,1) &   2  & Thm \ref{T=n-2}    \\
\hline
60& (0,\,1,\,1,\,1,\,0,\,1,\,1,\,1) &   2  & Thm \ref{T=n-2} \\
\hline
61& (0,\,1,\,1,\,0,\,1,\,1,\,1,\,1) &   2  & Thm \ref{T=n-2} \\
\hline
62& (0,\,1,\,0,\,1,\,1,\,1,\,1,\,1) &   2  & Thm \ref{T=n-2}    \\
\hline
63& (0,\,0,\,1,\,1,\,1,\,1,\,1,\,1) &   2  & Thm \ref{T=n-2} \\
\hline
64& (0,\,1,\,1,\,1,\,1,\,1,\,1,\,1) &   2  & Complete graph  \\
\hline
\end{tabular}}

\section{An upper bound for $q(G)$ of threshold graphs}\label{upper bound}

In this section, we provide a general upper bound for any connected threshold graph of order $n$.

\begin{theorem}
 Let $G\cong (\underbrace{0, \ldots, 0}_{k_1},\underbrace{1, \ldots, 1}_{t_1}, \underbrace{0, \ldots, 0}_{k_2},\underbrace{1, \ldots, 1}_{t_2}, \ldots,  \underbrace{0, \ldots, 0}_{k_s}, \underbrace{1, \ldots, 1}_{t_s})$ with $k_i, t_i, s\ge 1$ be a connected threshold graph with $n$ vertices. Then
 $$q(G)\le \left\lfloor\frac{s+9}{2} \right\rfloor,$$ which implies 
 $$q(G)\le \frac{n}{4}+C$$ for some constant $C$ independent of $n$.
\end{theorem} 
\begin{proof}
From Lemma \ref{jdup}, we have $q(G)\le q(H)$, where $H\cong (\underbrace{0, \ldots, 0}_{k_1},1, \underbrace{0, \ldots, 0}_{k_2},1, \ldots,  \underbrace{0, \ldots, 0}_{k_s}, 1)$. For the graph  $H$, if for any $i\in [s]$, $k_i=1$ then by Lemma \ref{thm00} $$q(H)=3\le \left\lfloor\frac{s+9}{2} \right\rfloor.$$ Otherwise, there is at least one index $i\in [s]$ with $k_i>1$. Consider the related threshold graph $H_0=(\underbrace{0,1,\ldots, 0,1}_{2p_1}, 0,0, 1,\underbrace{0,1,\ldots, 0,1}_{2p_2})$, where $p_1+p_2+1=s$. We have the following three cases:\\
 
 $Case\,1)\,$ Both $p_1$ and $p_2$ are even. Let $p_1=2s_1$, and $p_2=2s_2$. In this case, consider the following vertex-clique incidence matrix $M$:
 $$M=\left(
  \begin{array}{ccccccccccccccc}
    1 & 0 & 0 & 0 & 0 & 0 &\ldots &  & 0 & 0 \\
    1 & 2 & 0 & 0 &  0 & 0 &\ldots &  & 0 & 0 \\
    1 & 1 & 1 & 0 &  0 & 0 &\ldots &  & 0 & 0 \\
    1 & 1 & 1 & 2 & 0 & 0 & \ldots &  & 0 & 0 \\
    1 & 1 & 1 & 1 & 1 & 1 &\ldots &  & 0 & 0 \\
    1 & 1 & 1 & 1 &  1 & 1 &\ldots &  & 0 & 0 \\
    1 & 1 & \vdots & \vdots & \vdots & \vdots & \ldots &  & 0 & 0 \\
\vdots & \vdots & 1 & 1 &   &  &\ldots &  & 0 & 0 \\
    1 & 1 & 1 & 1 &  1 & 1 &\ldots &  & 1 & 0\\
    1 & 1 & 1 & 1 &  1 & 1 &\ldots &  & 1 & 2\\
\sqrt{3} & 0 & 0 & 0&  0 & 0 &\ldots &  & 0 & 0 \\
    0 & 1 & 0 & 0  & 0 & 0 & \ldots &  & 0 & 0 \\
    0 & 0 & \sqrt{3} & 0 & 0 & 0 & \ddots &  & 0 & 0 \\
    0 & 0 & 0 & 1 & 0 & 0 &  &  & 0 & 0 \\
    \vdots & \vdots  & \vdots  & \ddots & 1 & 0 &\ddots   &  & 0 & 0 \\
    0 & 0 & 0 &  \ldots & 0 & 1 &  &  & 0 & 0 \\
    0 & 0 & 0 &  \ldots & 0 & 0 & \ddots &  & 0 & 0 \\
    0 & 0 & 0 &  \ldots & 0 & 0 &  &  & \sqrt{3} & 0 \\
    0 & 0 & 0 &  \ldots & 0 & 0 & \ldots &  & 0 & 1 \\
  \end{array}
\right).$$
 
 Then
 $$B=M^T\,M=\left(
  \begin{array}{ccccccc}
    B_{11} & B_{12} &  &  & \ldots &  &  B_{1(s_1+s_2+1)} \\
    B_{12}^T & \ddots &   & \ddots  &  &  &\\
     &  &  B_{s_1s_1}  &  &    &  & \\
     &  &    & B_{(s_1+1)(s_1+1)} &  &  &  \vdots \\
    \vdots &  & \ddots  &  & B_{(s_1+2)(s_1+2)} & & \\
      &  &    &  &   &  \ddots &  \\
    B_{1(s_1+s_2+1)}^T  &  &  & \ldots &    & & B_{(s_1+s_2+1)(s_1+s_2+1)} \\
  \end{array}
\right),$$ where $B_{11}=\left( \begin{array}{cc}
     T+3 & T  \\
     T & T+3 \\
  \end{array}
\right)$, \, $B_{s_1s_1}=\left( \begin{array}{cc}
     p_2+6 & p_2+3  \\
     p_2+3 & p_2+6 \\
  \end{array}
\right)$, \, $B_{(s_1+1)(s_1+1)}=\left( \begin{array}{cc}
     p_2+2 & p_2+1  \\
     p_2+1 & p_2+2 \\
  \end{array}
\right)$, \, $B_{(s_1+2)(s_1+2)}=\left( \begin{array}{cc}
     p_2+1 & p_2-2  \\
     p_2-2 & p_2+1 \\
  \end{array}
\right)$, \, $B_{(s_1+s_2+1)(s_1+s_2+1)}=\left( \begin{array}{cc}
     5 & 2  \\
     2 & 5 \\
  \end{array}
\right)$ and $B_{ij}$ is a $2\times 2$ matrix with all entries equal if $i\neq j$. We have $\{3^{[{\frac{s-1}{2}}]},\, 1\}\subseteq \Spec(B),$ and then $\{0^{[s]},\, 3^{[{\frac{s-1}{2}}]},\, 1\}\subseteq \Spec(A)$, where $A=M\,M^T$. The total number of eigenvalues of $A$ is $2s+1$, so the number of remaining unaccounted for eigenvalues is $(s+1)/2$.
Thus  $$q(A)\le \frac{s+1}{2}+3=\frac{s+7}{2}.$$  Obviously $A=M\,M^T\in S(H_0)$ and the eigenvalues of $A$ (which are also eigenvalues of $B$) equal to $3$ and $1$ appear on the main diagonal of $A$. Then applying Lemma \ref{dup} we may duplicate vertices corresponding to the isolated vertices in $H_0$ and not increase the total number of  eigenvalues. This gives the threshold graph $H$ has a realization with the same set of eigenvalues as $H_0$. Hence,
$$q(G)\le q(H) \le \frac{s+7}{2}.$$

 $Case\,2)\,$ Both $p_1$ and $p_2$ are odd. Applying a similar method as above, we can show $q(G)\le \frac{s+9}{2}$.\\

 $Case\,3)\,$  $p_1$ and $p_2$ have different parity. Using similar techniques, it follows that  $q(G)\le \frac{s+8}{2}$.\\
 
 Hence the first part of the proof is done. We have $n\ge 2s$, i.e., $s\le \frac{n}{2}$. This gives $$q(G)\le \frac{s+9}{2} \le \frac{n+18}{4}=\frac{n}{4}+\frac{18}{4},$$ which completes the proof of the theorem. 
\end{proof}

\begin{theorem}
 Let $G$ be a graph of order $n$ with a cycle of length four and all other vertices adjacent to a vertex of the cycle. Then $q(G)=3$.
\end{theorem}
\begin{proof} Let $M$ be a vertex-clique incidence matrix of $G$ as follows:
$$M=\left(
  \begin{array}{cccccc}
    \sqrt{n-3} & 0 & 0 & 0 & \ldots & 0 \\
    -1 & 1 & 0 & 0 & \ldots & 0 \\
    1 & 1 & 0 & 0 & \ldots & 0 \\
    0 & 1 & 1 & 1 & \ldots & 1\\
    0 & 0 & \sqrt{2} & 0& \ldots & 0 \\
    0 & 0 & 0 & \sqrt{2} &  & \vdots \\
    \vdots & \vdots  & \vdots  & \ddots &\ddots   & 0 \\
    0 & 0 & 0 &  \ldots & 0 & \sqrt{2} \\
  \end{array}
\right).$$ We have $A=M M^T\in S(G)$. To find its spectrum, we determine the eigenvalues of the $(n-2)\times (n-2)$ matrix 
$$M^T M=\left(
  \begin{array}{cccccc}
    n-1 & 0 & 0 & \ldots & 0 \\
    0 & 3 & 1 & \ldots & 1 \\
    0 & 1 & 3 & \ddots & \vdots\\
    \vdots & \vdots & \ddots & \ddots & 1\\
    0 & 1 & \dots& 1 & 3 \\
  \end{array}
\right),$$ which are as follows: $$\Spec(M^TM)=\{n-1^{[2]},\,2^{[n-4]}\}.$$ Then the spectrum of $A$ is 
$$\Spec(A)=\{n-1^{[2]},\,2^{[n-4]},\, 0^2\}.$$ This with Lemma \ref{diam} gives $q(G)=3$. 
\end{proof}

\section{Conclusions}
This paper has laid some groundwork for studies on the spectral parameter $q(G)$ for certain structured families of graphs. The specific family considered here was threshold graphs, where we established numerous results concerning $q(G)$, including: threshold graphs on $n$ vertices with specific values of the trace; a general upper bound on $q$ for all threshold graphs; a workable constraint for threshold graphs with $q$ equal to two; and a complete list of $q(G)$ for threshold graphs of order 7 and 8 (for the latter there are just two exceptions). 

Looking ahead there are plenty of unresolved issues on this topic including:

\begin{enumerate}
\item Consider $q(G)$ for general threshold graphs with fixed trace $T$ with $T \in \{3,4,\ldots, n-4\}$;
\item Characterize the threshold graphs $G$ that satisfy $q(G)=3$ or $q(G)=4$, etc.;
\item Determine the maximum value of $q(G)$ for fixed $n$ and trace; and
\item Expand this study beyond the family of threshold graphs to include perhaps split graphs or cographs.
\end{enumerate}

\section*{Acknowledgements}
Dr. Fallat's research was supported in part by an NSERC Discovery Research Grant, Application No.: RGPIN-2019-03934.
The work of the PIMS Postdoctoral Fellow Dr. S. A. Mojallal leading to this publication was supported in part by PIMS. \\

\end{document}